\documentclass[a4paper, 11pt, english]{article}
\usepackage[utf8]{inputenc}
\usepackage[T1]{fontenc}
\usepackage{babel}
\usepackage{graphicx}
\usepackage{stmaryrd}
\usepackage[a4paper]{geometry}
\usepackage{amsmath,amsfonts,amssymb,amsthm,epsfig,epstopdf,url,array}
\usepackage{rotating}
\usepackage[colorlinks=true,citecolor=red,linkcolor=blue,pdfpagetransition=Blinds]{hyperref}
\usepackage{cleveref}
\usepackage{nameref}
\usepackage{enumitem}
\usepackage{comment}
\Crefname{paragraph}{Section}{Sections}
\usepackage{fancyhdr}
\usepackage{mathrsfs}

\usepackage{fullpage}

\newcommand{\ensemblenombre}[1]{\mathbb{#1}}
\newcommand{\N}{\ensemblenombre{N}}
\newcommand{\Z}{\ensemblenombre{Z}}

\newcommand{\R}{} 
\renewcommand{\R}{\ensemblenombre{R}}

\newcommand{\T}{\ensemblenombre{T}}

\renewcommand{\leq}{\leqslant}
\renewcommand{\geq}{\geqslant}

\newcommand{\dive}[1]{\mathrm{div}}

\newcommand{\eps}{\varepsilon}

\theoremstyle{plain} 
\newtheorem{pr}{Proposition}[section] 
\newtheorem{tm}[pr]{Theorem}
\newtheorem{lm}[pr]{Lemma}
\newtheorem{cor}[pr]{Corollary}
\theoremstyle{definition}
\newtheorem{df}[pr]{Definition}
\newtheorem{rmk}[pr]{Remark}

\numberwithin{equation}{section}

\allowdisplaybreaks

\makeatletter
\let\original@addcontentsline\addcontentsline
\newcommand{\dummy@addcontentsline}[3]{}
\newcommand{\DeactivateToc}{\let\addcontentsline\dummy@addcontentsline}
\newcommand{\ActivateToc}{\let\addcontentsline\original@addcontentsline}
\makeatother

\begin{document}

	\title{Quantitative observability of the Schrödinger equation in the torus from open sets}
	\author{Kévin Le Balc'h, Jiaqi Yu}
	
	\maketitle
	
	\begin{abstract}
		We study the observability of the Schrödinger equation on the $d$-dimensional torus $\T^d=\R^d/(2\pi\Z)^d$, $d \geq 1$, from an open subset $\omega \subset \T^d$. Our main result establishes a quantitative observability estimate for the free Schrödinger equation in the regime of small times $T$ and small observation sets of the form $\omega = \prod_{j=1}^{d}(a_j,b_j)$. A central ingredient in the proof is a cluster decomposition method combined with an induction scheme introduced by Bourgain and further developed by Burq and Zhu.
	\end{abstract}

	
	\section{Introduction}

	The observability of the Schrödinger equation on the $d$-dimensional torus $\T^d=\R^d/(2\pi\Z)^d$ from an open subset $\omega \subset \T^d$ is a central topic in control theory and mathematical physics. It concerns the ability to recover global information about the solution by monitoring it only on a subregion over a finite time interval. 
	
	\medskip
	
	In this paper, we consider the free Schr\"{o}dinger equation
	\begin{align}
		\label{eq:schro}
		\left\{
		\begin{array}{ll}
			i\partial_t u - \Delta u = 0, & (t,x)\in \R \times \T^d, \\
			u(0,\cdot) = u_0 \in L^2(\T^d).
		\end{array}
		\right.
	\end{align}

	Our main result provides a quantitative observability estimate for \eqref{eq:schro} when the observation set $\omega \subset \T^d$ is a product of (small) intervals and the time interval $(0,T)$ is small.
	
	\begin{tm}\label{Thm1}
		There exist constants $C,C'>1$, depending only on $d$, such that for every $\omega=\prod_{j=1}^{d}(a_j,b_j)\subset\T^d$ with $0<|\omega|<1$, for every $0<T<1$ and for every $u_0 \in L^2(\T^d)$,
		\begin{equation}\label{equ-obs}
			\|u_0\|_{L^2(\T^d)}^2 \leq C_{\text{obs}}(T,\omega)
			\int_0^T \|e^{-it\Delta}u_0\|_{L^2(\omega)}^2 \, dt,
		\end{equation}
		where
		\begin{equation}
			\label{eq:defCTomega}
			C_{\mathrm{obs}}(T,\omega)
			=
			C\underbrace{\exp\left\{
				\exp\left\{
				\cdots
				\exp\left\{
				\frac{C}{(T|\omega|)^{C'}}
				\right\}
				\cdots
				\right\}
				\right\}}_{d\text{ times}}.
		\end{equation}
	\end{tm}
	
	\begin{rmk}
		Of course, \Cref{Thm1} also applies to any (small) ball $\omega = B(x_0,r) \subset \T^d$ by choosing a rectangle contained in that ball, leading to 
		\begin{equation}\label{equ-obsball}
			\|u_0\|_{L^2(\T^d)}^2 \leq C_{\text{obs}}(T,r)
			\int_0^T \|e^{-it\Delta}u_0\|_{L^2(\omega)}^2 \, dt,
		\end{equation}
		where
		\begin{equation}
			\label{eq:defCTr}
			C_{\mathrm{obs}}(T,r)
			=
			C\underbrace{\exp\left\{
				\exp\left\{
				\cdots
				\exp\left\{
				\frac{C}{(Tr)^{C'}}
				\right\}
				\cdots
				\right\}
				\right\}}_{d\text{ times}}.
		\end{equation}
		But we choose to keep the formulation with $\omega$ a product of small intervals, because as we will see, this particular form of observation sets play a key role in the proof based on induction on the dimension $d$.
	\end{rmk}

	\paragraph*{Bibliographical comments on the observability of the Schrödinger equation in the torus.}
	
	The observability of the Schrödinger equation on the torus $\T^d$ from arbitrary open observation sets $\omega$ has a long story that goes back at the early nineties. We briefly review some known results.
	
	\medskip
	
	Jaffard \cite{Jaf90} established the observability estimate \eqref{equ-obs} for the free Schrödinger equation \eqref{eq:schro} in dimension $d=2$, for arbitrary nonempty open sets $\omega$. His proof relies on lacunary Fourier series techniques developed by Kahane \cite{Kah62}. Komornik and Loreti \cite{KL05} later extended this result to all dimensions $d \geq 1$, still in the free case \eqref{eq:schro}.
	
	Burq and Zworski \cite{BZ12} proved observability for the Schrödinger equation with real-valued time independent potentials $V=V(x)$
	\begin{align}
		\label{eq:schroV}
		\left\{
		\begin{array}{ll}
			i\partial_t u - \Delta u + V u = 0, & (t,x)\in \R \times \T^d, \\
			u(0,\cdot) = u_0 \in L^2(\T^d).
		\end{array}
		\right.
	\end{align}
	in dimension $d=2$, for potentials $V \in C^0(\T^2)$ and arbitrary nonempty open sets $\omega$. Their approach is based on semiclassical analysis combined with a normal form argument. Around the same time, Anantharaman and Macià \cite{AM14} established observability for \eqref{eq:schroV} in all dimensions $d \geq 1$, for arbitrary open sets $\omega$ and bounded Riemann-integrable potentials $V$, using semiclassical measures and quantum limits.
	
	Bourgain, Burq, and Zworski \cite{BBZ13} subsequently addressed the case of $L^2(\T^2)$ potentials in dimension $d=2$, with Strichartz estimates playing a key role in their analysis. Bourgain \cite{Bou14} further investigated the three-dimensional case $d=3$ with $L^\infty(\T^3)$ potentials, introducing a novel approach based on cluster decomposition techniques. Tao \cite{Tao21} obtained quantitative observability estimates for the free equation \eqref{eq:schro} in dimension $d=2$, for small balls $\omega = B(x_0, \varepsilon)$ and time $T = 2\pi$, refining the strategy developed in \cite{BBZ13}.
	
	Burq and Zworski \cite{BZ19} returned to the free case $V=0$ in dimension $d=2$, establishing observability estimates for arbitrary measurable sets $\omega$ with $|\omega| > 0$. Their proof combines semiclassical analysis with Strichartz estimates to handle rough observation sets. The first author and Martin \cite{LBM23} extended these results to potentials $V \in L^\infty(\T^2)$ in dimension $d=2$, again for observation sets of positive measure.
	
	More recently, Burq and Zhu \cite{BZ25} investigated observability for rough sets and a broad class of potentials, including Riemann-integrable ones, using an approach inspired by Bourgain’s cluster decomposition method. Let us also mention the recent work of Burq, Germain, Sorella and Zhu \cite{BGSZ26}, who study trace and observability inequalities for Laplace eigenfunctions on $\mathbb T^d$ with respect to general Borel measures. Their results provide a complementary stationary and spectral viewpoint on observability on the torus.
	
	In dimension $d=1$, Alphonse and Tzvetkov \cite{AT25} revisited observability from rough sets, notably establishing a smoothing effect for the fractional Schrödinger equation, which is of independent interest. Alphonse and Lafontaine \cite{AL25} proved observability estimates for the Schrödinger equation on the equilateral triangle, under both Neumann and Dirichlet boundary conditions and for open observation sets. Their strategy relies on Pinsky's tiling argument, reducing the problem to observability estimates on rational twisted tori.
	
	It is worth emphasizing that, despite the diversity of techniques developed to study observability for the Schrödinger equation, many of these works rely, either explicitly or implicitly, on a dimensional reduction argument introduced in this context.
	
	Beyond the flat torus, the validity of observability from arbitrary nonempty open sets is strongly influenced by the underlying geometry. In the setting of compact surfaces with Anosov geodesic flow, Dyatlov, Jin and Nonnenmacher \cite{DJN22} established uniform lower bounds on the mass of high-frequency eigenfunctions on every nonempty open set, which in particular imply controllability for the Schrödinger equation.
	
	\paragraph*{Comparison of our main result with the existing literature.}
	
	Our main result from \Cref{Thm1} is the first result in the literature establishing a quantitative observability estimate of the free Schrödinger equation \eqref{eq:schro} in the multi-dimensional torus in the regime of small rectangle observation sets and small times.
	
	In dimension $1$, our estimate takes the following form
	\begin{equation}\label{equ-obs1d1}
		\|u_0\|_{L^2(\T)}^2 \leq C \exp\left( \frac{C}{(|\omega| T)^{C'}}\right)
		\int_0^T \|e^{-it\Delta}u_0\|_{L^2(\omega)}^2 \, dt.
	\end{equation}
	This estimate is far from being sharp according to \cite[Theorems 3.1 and 2.2]{Mil05} where the author is able to prove	
	\begin{equation}\label{equ-obs1d1Miller}
		\|u_0\|_{L^2(\T)}^2 \leq \exp\left( \frac{32 \pi^2}{T}\right)
		\int_0^T \|e^{-it\Delta}u_0\|_{L^2(\omega)}^2 \, dt.
	\end{equation}
	This estimate \eqref{equ-obs1d1Miller} is actually very specific to the one-dimensional case because it mainly relies on the transmutation method that consists in deducing an observability estimate for the Schrödinger equation from the wave equation. Note also that the asymptotic dependence of \eqref{equ-obs1d1Miller} in function of small $T$, up to universal positive constants (see also \cite{TT07, TT09, Lis14}) is known to be sharp in such a geometric configuration \cite{Mil04}. Last but not least, it is interesting to note that the observability constant \eqref{equ-obs1d1Miller} does not blow-up as $|\omega| \to 0$, contrary to our estimate \eqref{equ-obs1d1}. Again this appears to be an artefact of the one-dimensional case for which the geometric control condition is always satisfied.

	In dimension $2$, recall that Tao \cite{Tao21} considered the case of arbitrary small ball $\omega = B(x_0, \varepsilon)$ and obtained the observability estimate at time $T=2 \pi$,
	\begin{equation}\label{equ-obsTao}
		\|u_0\|_{L^2(\T^2)}^2 \leq C \exp\left(\exp\!\left(\frac{C\log \eps^{-1}}{\log\log\varepsilon ^{-1}}\right)\right)
		\int_0^{2 \pi} \|e^{-it\Delta}u_0\|_{L^2(B(x_0, \varepsilon))}^2 \, dt.
	\end{equation}
	In this case, our observability estimate \eqref{equ-obsball} leads to
	\begin{equation}\label{equ-obs1}
		\|u_0\|_{L^2(\T^d)}^2 \leq C \exp\left( \exp\!\left(\frac{C}{(\varepsilon T)^{C'}}\right)\right)
		\int_0^T \|e^{-it\Delta}u_0\|_{L^2(\omega)}^2 \, dt.
	\end{equation}
	While his estimate \eqref{equ-obsTao} is asymptotically better than ours \eqref{equ-obs1} as $\varepsilon \to 0$, we notice that in such a case, both observability constants blow-up as $|\omega| \to 0$. As previously remarked by Tao, the behaviour of the constant in function of $T$ and $\omega$ mainly comes from the low frequencies parts. Despite its importance, this problematic has not been too much considered in the literature for the Schrödinger equation (see the recent article \cite{SSY25} in the real line), in comparison to what is known for the heat equation, see for instance \cite{LL21} and the references therein.
	

	\paragraph*{Applications of our main result.} Our observability estimate coming from \Cref{Thm1} has standard applications as resolvent estimates for the (stationary) Schrödinger operator \cite{BZ04}, \cite{Mil05} and controllability results for the (controlled) Schrödinger equation
	\begin{align}
		\label{eq:schroCont}
		\left\{
		\begin{array}{ll}
			i\partial_t \psi - \Delta \psi = h 1_{\omega}, & (t,x)\in \R \times \T^d, \\
			\psi(0,\cdot) = \psi_0 \in L^2(\T^d).
		\end{array}
		\right.
	\end{align}
	by the Hilbert Uniqueness Method, see for instance \cite{Cor07}. We only state the exact controllability result explicitly.

	\begin{cor}
	There exist constants $C,C'>1$, depending only on $d$, such that for every $\omega=\prod_{j=1}^{d}(a_j,b_j)\subset\T^d$ with $0<|\omega|<1$, for every $0<T<1$ and for every $\psi_0, \psi_1 \in L^2(\T^d)$, there exists $h \in L^2(0,T;L^2(\omega))$ satisfying
	\begin{equation}
		\|h\|_{L^2(0,T;L^2(\omega))} \leq C_{\text{obs}}(T,\omega) (\|\psi_0\|_{L^2(\T^d)} + \|\psi_1\|_{L^2(\T^d)}),
	\end{equation}
	where $C_{\text{obs}}(T,\omega)$ is as in \eqref{eq:defCTomega} and such that the solution $\psi$ of \eqref{eq:schroCont} satisfies $$\psi(T,\cdot) = \psi_1.$$
	\end{cor}

	\paragraph*{Strategy of the proof.} Our main tools for proving \Cref{Thm1} are the cluster decomposition technique and a mathematical induction argument, borrowed from the papers \cite{Bou14} and \cite{BZ25}.
	
	We first reduce \Cref{Thm1} to \Cref{thm:obslowdim}, see below.
	The proof of \Cref{thm:obslowdim} is divided into two parts.
	
	The first part consists in proving a weak observability estimate for initial data with appropriate location of the Fourier support
	(\Cref{prop-weak-obs}). Following the strategy of \cite{BZ25}, we
	exploit the cluster structure and derive the estimate by induction on the dimension. A new observation is that, thanks to the Galilean
	transform (\Cref{prop:equiv}) and the analysis of the lattice structure (\Cref{lem:lattice}), one may reduce to the case where the remainder term is supported in a low-frequency region (\Cref{lem:support}) .
	
	The second part consists in absorbing the remainder term by a quantitative unique continuation estimate (\Cref{prop:lowfreq}). For this, we follow the framework of \cite[Proposition 5.1]{SSY25} and reduce the problem to an observability estimate for the heat equation (see \Cref{prop:Heat}), established via a quantitative Lebeau--Robbiano argument.

	\paragraph*{Acknowledgements.} Both authors would like to warmly thank Chenmin Sun for suggesting the study of quantitative observability estimates of the Schrödinger equation in the torus and many discussions, comments around this article. Both authors would like to also thank Nicolas Burq and Hui Zhu for addressing some comments about a first version of this article.
	
	\section{Preliminary tools}
	
	For any $f\in L^2(\T^d)$, $\Gamma\subset\Z^d$, we define the Fourier projector
	\[
	\Pi_\Gamma f=\sum_{n\in\Gamma} \hat{f}(n)e^{in\cdot x},
	\]
	where
	\[
	\hat{f}(n)=(2 \pi)^{-d/2} \int_{\T^d}f(x)e^{-in\cdot x}dx,\quad n\in\Z^d.
	\]
	Let $L_\Gamma^2=\Pi_\Gamma(L_x^2)$. In the following, for any $S\subset\Z^d$, we denote by dim $S$ the dimension of the smallest affine sublattice in $\Z^d$ that contains $S$.
	
	\medskip
	
	\subsection{From open sets to a band-limited smooth function}\label{sect:psi}
	
	We begin with the following approximation's result on the one-dimensional torus of the indicator function $\chi_{I_j}$ where $I_j$ is some interval contained in $\T$.
	
	\begin{lm}\label{lem-psi}
		Let $I_j\subset\T$ be an interval. There exist $\varphi_j\in C^{\infty}(\T)$ such that 
		\begin{enumerate}[label=(\roman*)]
			\item $0 \leq \varphi_j \leq 1$,
			\item $\mathrm{supp}(\widehat{\varphi_j}) \subset (-R, R)$,
			\item $\|\varphi_j-\chi_{I_j}\|_{L^1(\T)}\leq C_0R^{-1/4}$ for all $R>4$ where $C_0>0$ is an universal constant.
		\end{enumerate}
	\end{lm}
	\begin{proof}
		We do a frequency regularization on the torus $\T$. 
		For fixed $R\in\mathbb N$, consider a Fejér kernel $K_R\in C^\infty(\mathbb T)$
		satisfying
		\[
		K_R\geq 0,\qquad \int_{\mathbb T} K_R(z)\,dz = 1,
		\qquad
		\widehat{K_R}(\xi)=0\ \text{ if }|\xi|>R,\quad \xi\in\mathbb Z.
		\]
		Now, we define
		\[
		\varphi_j(s) := (K_R * \chi_{I_j})(s), 
		\qquad s\in\mathbb T,
		\]
		the convolution being taken on the torus $\mathbb T$.
		
		The regularity property $\varphi_j\in C^{\infty}(\T)$ comes from convolution property and $K_R \in C^\infty(\mathbb T)$.
		
		Since $K_R\geq 0$, $\int K_R=1$ and $0\leq \chi_{I_j} \leq 1$ we have that $(i)$ is satisfied.  
		
		The support property $(ii)$ comes from $\widehat{\varphi_j}(\xi) = \widehat{K_R}(\xi)\,\widehat{\chi_{I_j}}(\xi)$, hence
		$\operatorname{supp}\widehat{\varphi_j} 
		\subset \operatorname{supp}\widehat{K_R} \subset (-R,R)$.

		We now assume $I_j=(a_j,b_j)$ with $-\pi<a_j<b_j<\pi$. Choosing $\phi_{j,n}(s)\in C_0^\infty(\T)$ such that $\phi_{j,n}=0$ on $\T\setminus (a_j-\frac{1}{n},b_j+\frac{1}{n})$ and $\phi_{j,n}=1$ on $(a_j+\frac{1}{n},b_j-\frac{1}{n})$ with $|\phi_{j,n}^{(k)}(s)|\leq C_kn^k$ for all $k \in \N$. We then have
		\[
		\|\chi_{I_j}-\phi_{j,n}\|_{L^1(\T)}\leq \frac{4}{n},
		\]
		and by \cite[p.162 (Theorem 2)]{Na64}
		\[
		\|K_R*\phi_{j,n}-\phi_{j,n}\|_{L^\infty(\T)}<\frac{C n\ln R}{R},\qquad R>1.
		\]Thus,
		\begin{align*}
			\|K_R*\chi_{I_j}-\chi_{I_j}\|_{L^1(\T)}&\leq (\sup_R\|K_R\|_{L^1(\T)}+1)\|\chi_{I_j}-\phi_{j,n}\|_{L^1(\T)}+\|K_R*\phi_{j,n}-\phi_{j,n}\|_{L^\infty(\T)}\\
			&\leq \frac{C}{n}+\frac{Cn\ln R}{R}.
		\end{align*}
		Taking $n=R^{1/4}$, we obtain $(iii)$.
	\end{proof}		
	
	We then define
	$$ I_j = (a_j,b_j) \subset \T,\qquad j \in \{1,\dots, d\},$$
	and assuming that $0 < T < 2 \pi$, one can identify $\T_t$ with $(0,2 \pi)$ so that the interval $(0,T)$ belongs to $\T_t$ and we set
	$$ I_{d+1} = (0,T) \subset \T_t.$$
	From the previous lemma, we can now define 
	\begin{equation}
		\label{eq:defpsi}
		\psi(x,t)=\big(\prod_j^d\varphi_j(x_j)\big) \varphi_{d+1}(t)\qquad (x,t) \in \T^d \times \T_t.
	\end{equation}
	Notice that $\widehat\psi(\xi,\tau)=\big(\prod_j^d\widehat\varphi_j(\xi_j)\big) \widehat\varphi_{d+1}(\tau)$, where $(\xi,\tau)=(\xi_1,\cdots,\xi_d,\tau)\in\Z^{d+1}$. So,
	\begin{equation}
		\label{eq:supppsi}
		\mathrm{supp}\ \widehat{\psi}\subset(-R,R)^{d+1}.
	\end{equation}

	\subsection{Lattice properties}
	We recall the following concepts (see \cite{Cas97}):
	\begin{df}\label{def:lattice}
		Let $\textbf{u}^1,\,\textbf{u}^2,\cdots,\textbf{u}^d$ be $d$ independent points in $\R^d$. Then the set $\Lambda$ of points
		\[
		x=a_1\textbf{u}^1+\cdots+a_d\textbf{u}^d\qquad\qquad(a_1,\cdots,a_d\quad \mbox{ integers})
		\]
		is called a lattice. The system of points $\textbf{u}^1,\cdots,\textbf{u}^d$ is called a basis of $\Lambda$.
	\end{df}
	
	\begin{df}\label{def:primi}
		A system of $k\leq n$ independent points $\textbf{v}^1,\cdots,\textbf{v}^k$ in a lattice $\Lambda$ is called primitive, if the points of $\Lambda\cap \text{span}_\R(\textbf{v}^1,\cdots,\textbf{v}^k)$ are all of the form:
		\[
		b_1\textbf{v}^1+\cdots+b_k\textbf{v}^k\qquad\qquad(b_1,\cdots,b_k\quad \mbox{ integers}).
		\]
	\end{df}
	
	Let $\textbf{V}=(\textbf{v}^1,\cdots,\textbf{v}^k)\in\R^{d\times k}$, and
	\begin{equation*}
		\mathcal{P}(\textbf{V}):=\Big\{\sum_{j=1}^{k}t_j\textbf{v}^j,\quad0\leq t_j<1\Big\}.
	\end{equation*}
	\begin{df}\label{def:det}
		Let $\Lambda=\Lambda(\textbf{V})\subset\R^d$ for some basis $\textbf{V}\in \R^{d\times k}$. Define the determinant of $\Lambda$ to be $\det(\Lambda)=\sqrt{\det(\textbf{V}^{\mathsf{T}} \textbf{V})}$.
	\end{df}
	
	We need the standard notion of the \emph{successive minima} of a lattice $\Lambda$ in $\R^d$. For $j=1,...,d$, we recall that $\mu_j(\Lambda)$ is the smallest real number $\mu$ such that $\Lambda$ contains $j$ linearly independent vectors of length at most $\mu$.
	
	Now, suppose that $\Lambda\subset\Z^d$ is a primitive lattice, i.e., $\text{span}_{\R}(\Lambda)\cap\Z^d=\Lambda$ where $\text{span}_{\R}(\Lambda)$ is the
	linear span of $\Lambda$ in $\R^d$. Let $\Gamma$ be the affine sublattice in $\mathbb{A}_\Lambda:=\{q+\Lambda:\,q\in\Z^d\}$. Let $\Lambda^\perp:=\{k\in\Z^d:\,k\perp \Lambda\}$. Denote the orthogonal projection on $\text{span}_\R(\Lambda^\perp)$ by $\text{proj}_{\Lambda^\perp}$.
	
	\begin{lm}\label{lem:lattice}
		Suppose that $\Gamma=q+\Lambda\in\mathbb{A}_\Lambda$, and $\Lambda$ is a $s$-dimensional primitive lattice with $s\leq d$.  
		Then, there is $n_\Gamma\in\Lambda^\perp$ such that
		\begin{equation*}
			\mathrm{dist}(q+n_\Gamma,\,\mathrm{span}_\R(\Lambda))\leq C\det(\Lambda),
		\end{equation*}
		where $C$ depends only on $d,s$.
	\end{lm}
	
	\begin{proof}
		If $s=d$, then $\Lambda^\perp=\{0\}$, and
		\[
		\text{dist}(q,\,\text{span}_\R(\Lambda))=0.
		\] 
		Thus, the conclusion holds trivially when $s=d$.
		
		Now, we suppose $s\leq d-1$.
		Notice that $\Lambda^\perp$ is a $(d-s)$-dimensional sublattice of $\Z^d$. Let $r:=d-s\geq 1$.
		By the definition of the successive minima, there exist linearly independent vectors $\textbf{v}^1,\cdots,\textbf{v}^r\in\Lambda^\perp$ such that
		\begin{equation*}
			|\textbf{v}^j|\leq \mu_j(\Lambda^\perp),\qquad j=1,\cdots,r.
		\end{equation*}
		It follows from \cite[Lemma 2.1 and (4.3)]{DSS24} that 
		\begin{equation}\label{eq:minkow}
			|\textbf{v}^1|\cdots|\textbf{v}^r|\leq \mu_1(\Lambda^\perp)\cdots\mu_r(\Lambda^\perp)\leq c_{d,s}\det(\Lambda^\perp)^{r/d}.
		\end{equation}

		Let $N_0=\det(\Lambda)^2$. Consider now the group action 
		$$ (p, q + \Lambda) \in \Lambda^\perp \times \mathbb{A}_{\Lambda} \mapsto p + q + \Lambda.$$
		It follows from \cite[Lemma 3.2]{BZ25} that this group action has a finite number of orbits and we have
		\[
		\sharp(\mathbb{A}_\Lambda/\Lambda^\perp)=N_0.
		\]
		Thus, there exists $j\in\{1,\cdots,N_0\}$ such that
		\begin{equation}\label{eq:equiv-class}
			q+\Lambda=p+ r_j+\Lambda, \quad\mbox{ where }p\in \Lambda^\perp,\,r_j\in\Z^d.
		\end{equation}
		
		Let $\textbf{V}^\perp=\{\textbf{v}^1,\cdots,\textbf{v}^r\}$.
		We claim that there exists $r_j\in\Z^d$ satisfying $\text{proj}_{\Lambda^\perp}r_j\in\mathcal{P}(\textbf{V}^\perp)$ and \eqref{eq:equiv-class}. 
		
		Indeed, 
		since $\text{span}_{\Z}(\textbf{V}^\perp)$ is a sublattice of $\Lambda^\perp$, we have
		\[
		\text{span}_\R(\Lambda^\perp)=\bigcup\limits_{z\in\Lambda^\perp}(\mathcal{P}(\textbf{V}^\perp)+z).
		\]
		Thus, for every $x\in\text{span}_\R(\Lambda^\perp)$, there exists $z\in \Lambda^\perp$ such that $x-z\in \mathcal{P}(\textbf{V}^\perp)$.
		Suppose that $\tilde{r}_j$ satisfies \eqref{eq:equiv-class}, by the fact 
		$\text{proj}_{\Lambda^\perp}\tilde{r}_j\in \text{span}_\R(\Lambda^\perp)$, there is $z_j\in\Lambda^\perp $ such that $\text{proj}_{\Lambda^\perp}\tilde{r}_j-z_j\in \mathcal{P}(\textbf{V}^\perp)$. Taking $r_j=\tilde{r}_j-z_j$, then
		\[
		\text{proj}_{\Lambda^\perp}{r}_j=\text{proj}_{\Lambda^\perp}\tilde{r}_j-z_j \in \mathcal{P}(\textbf{V}^\perp),
		\]
		and 
		\[
		q+\Lambda=p+ \tilde{r_j}+\Lambda=p+r_j+ z_j+\Lambda:=\tilde{p}+r_j+\Lambda,
		\]
		where $\tilde{p}\in \Lambda^\perp,\,r_j\in\Z^d.$ Hence, there exists $r_j$ satisfying $\text{proj}_{\Lambda^\perp}r_j\in\mathcal{P}(\textbf{V}^\perp)$ and \eqref{eq:equiv-class}.
		
		Now, we choose $n_\Gamma=-p\in\Lambda^\perp$ with $p$ given by \eqref{eq:equiv-class}. Then
		\begin{multline}\label{eq:upp-diam}
			\quad\text{dist}(q+n_{\Gamma},\,\text{span}_\R(\Lambda))=\text{dist}(q-p,\,\text{span}_\R(\Lambda))=\text{dist}(r_j,\,\text{span}_\R(\Lambda))=|\text{proj}_{\Lambda^\perp}r_j|\\
			\leq \text{diam}(\mathcal{P}(\textbf{V}^\perp)),
		\end{multline}
		where the last inequality follows from $\text{proj}_{\Lambda^\perp}r_j\in\mathcal{P}(\textbf{V}^\perp)$.
		On the other hand,
		\begin{equation}\label{eq:upp-det}
			\text{diam}(\mathcal{P}(\textbf{V}^\perp))\leq \sum_{j=1}^r|\textbf{v}^j|\leq \sum_{j=1}^r\mu_j(\Lambda^\perp)\leq r\mu_r(\Lambda^\perp)\leq (d-s)c_{d,s}\det(\Lambda^\perp)^{r/d}.
		\end{equation}
		In \eqref{eq:upp-det}, the third inequality follows from $\mu_1(\Lambda^\perp)\leq \mu_2(\Lambda^\perp)\leq \cdots\leq \mu_r(\Lambda^\perp)$, and the last inequality follows from \eqref{eq:minkow} and the fact $\mu_j(\Lambda^\perp)\geq1$ for any $j$ recalling $\Lambda^\perp \subset \Z^d$. 
		Since $\Lambda$ is primitive, we get $\det(\Lambda^\perp)=\det(\Lambda)$ (see \cite{McM84}).
		This, combining with \eqref{eq:upp-diam}, \eqref{eq:upp-det}, completes the proof.
	\end{proof}

	\subsection{Galilean invariance}\label{subsect:Galilean} 
	We recall Galilean invariance property from \cite{BZ25}. For $p\in\R^d$, the \emph{Galilean transform} $\mathcal{G}_p$ is defined by
	\[
	(\mathcal{G}_p u)(t,x)=e^{i(p\cdot x+|p|^2t)}u(t,x+2pt).
	\]
	If $p\in \Z^d$ and $u(t)=e^{-it\Delta}u_0$, then
	\begin{equation}\label{equ-Gal-equiv}
		e^{-it\Delta}(e^{ip\cdot x}u_0(x))=(\mathcal{G}_p u)(t,x).
	\end{equation}
	If $\Gamma\in\mathbb{A}_\Lambda,\,p\in\Lambda^\perp$, and $u\in L_{\Z\times\Gamma}^2$, it follows from \cite[Lemma 3.1]{BZ25} that
	\begin{equation}\label{equ-Gal-1}
		|\mathcal{G}_p u|=|u|.
	\end{equation}
	From \eqref{equ-Gal-equiv} and \eqref{equ-Gal-1} we have
	\begin{pr}\label{prop:equiv}
		Let $\Gamma\in\mathbb{A}_\Lambda$, $p\in\Lambda^\perp$. The following statements are equivalent:
		
		\noindent(i) For any $T>0$ and each $u_0\in L_\Gamma^2(\mathbb T^d)$,
		there exists a positive constant $C_{1}(T,\omega)$ so that 
		\begin{equation*}
			\|u_0\|_{L^2(\T^d)}^2 \leq C_{1}(T,\omega)
			\int_0^T \|e^{-it\Delta}u_0\|_{L^2(\omega)}^2 \, dt.
		\end{equation*}
		
		\noindent(ii) For any $T>0$ and each $u_0\in L_\Gamma^2(\mathbb T^d)$,
		there exists a positive constant $C_{2}(T,\omega)$ so that 
		\begin{equation*}
			\|e^{ip\cdot x}u_0\|_{L^2(\T^d)}^2 \leq C_{2}(T,\omega)
			\int_0^T \|e^{-it\Delta}(e^{ip\cdot x}u_0)\|_{L^2(\omega)}^2 \, dt.
		\end{equation*}
		
		Furthermore, when one of the above two statements holds, the constants $C_{1}(T,\omega)$ and $C_{2}(T,\omega)$ can be chosen as the same  number.
	\end{pr}

	\subsection{The cluster decompositions}\label{subsect:cluster}
	
	In this part, we introduce the cluster decomposition. 
	
	\medskip

	Let $\Gamma\in\mathbb{A}_\Lambda$ be an affine sublattice, consider the set
	\begin{equation}\label{def:S_N}
		\Sigma(\Gamma) := \{(|n|^2,n)\in\mathbb{Z}\times\mathbb{Z}^d : n\in\Gamma\}.
	\end{equation}
	By a cluster decomposition lemma, see \cite[Lemma 3.3]{BZ25}, there exists $C_1>0$ depending on the dimension $d$, such that for $R>1$, 
	\begin{itemize}
		\item $\Sigma(\Gamma)$ is partitioned into disjoint clusters $Q_\Gamma^\alpha$, i.e. $$\Sigma(\Gamma) = \bigsqcup_{\alpha\in\mathscr{I}(\Gamma;R)} Q_\Gamma^\alpha,$$
		\item the diameter of each $Q_\Gamma^\alpha$ is uniformly bounded by $R^{C_1}$, i.e.
		\begin{equation}\label{equ-Cluster-diam}
			\mathrm{diam}(Q_\Gamma^\alpha) \leq R^{C_1},
		\end{equation}
		\item clusters are at least separated by a distance $100R$, i.e.
		\begin{equation}\label{equ-Cluster-dist}
			\mbox{dist}(Q_\Gamma^\alpha,Q_\Gamma^\beta)>100R, \mbox{ if }\alpha\neq\beta,
		\end{equation}
	\end{itemize}
	Let $Z_\Gamma^\alpha := \{ n\in\mathbb{Z}^d : (|n|^2,n)\in Q_\Gamma^\alpha\}$ and $\Gamma^\alpha$ be the smallest affine sublattice containing $Z_\Gamma^\alpha$. Denote the orthogonal projection on $\text{span}_\R(\Lambda)$ by $\text{proj}_\Lambda$. Let $\mathscr{I}_\flat(\Gamma;R)$ resp. $\mathscr{I}_\sharp(\Gamma;R)$ be the set of all $\alpha\in \mathscr{I}(\Gamma;R)$ such that $\text{dim } \Gamma^\alpha=\text{dim } \Gamma$ resp. $\text{dim } \Gamma^\alpha<\text{dim } \Gamma$. 
	Assume that $\Gamma\in \mathbb{A}_\Lambda$ and $\text{dim } \Gamma=s\geq1$. For every $u_0\in L_\Gamma^2$, set
	$u_{0,\alpha}=\sum_{n\in Z_\Gamma^\alpha} c_n e^{i n\cdot x}$.

	\medskip
	
	We have the following lemma.
	\begin{lm}
	\label{lem:support}
Let $u_0\in L^2_{\Gamma}$, where $\Gamma=n_\Gamma+q+\Lambda$ with $n_\Gamma\in \Lambda^\perp$ be given by \Cref{lem:lattice}. Then,
		\begin{align}
		\label{eq:estimationsupport}
			\text{supp }\hat{u}_{0,\alpha}\subset\{n\in\Z^d:\,|n|\leq c_1R^{c_2}\},\qquad \forall \alpha\in \mathscr{I}_\flat(\Gamma;R),
		\end{align}
		where $c_1,c_2>0$ depend only on $s,d$. 
	\end{lm}

	\begin{proof}
	By \cite[Lemma 3.4]{BZ25}, there exist $C_2,C_3>0$ depending only on $d$ such that, if $\alpha\in\mathscr{I}_\flat(\Gamma;R)$, then
	\begin{equation}
		\label{eq:estimateprojection}
		\sup_{n\in Z_\Gamma^\alpha}|\text{proj}_\Lambda\, n|<C_2R^{C_3}.
	\end{equation}
	Thus, for any $\alpha\in\mathscr{I}_\flat(\Gamma;R)$, $n\in Z_\Gamma^\alpha=\text{supp }\hat{u}_{0,\alpha}$, we have
		\begin{align}\label{equ-w2}
			|n|&\leq |\text{proj}_\Lambda\, n|+|\text{proj}_{\Lambda^\perp}\, n|=|\text{proj}_\Lambda\, n|+ \text{dist}(n,\text{span}_\R(\Lambda))\notag\\
			&\leq C_2R^{C_3}+\text{dist}(n,\text{span}_\R(\Lambda))\qquad\quad\mbox{(by \eqref{eq:estimateprojection})}\notag\\
			&=C_2R^{C_3}+\text{dist}(n_\Gamma+q,\text{span}_\R(\Lambda))\qquad\quad\mbox{(since $n\in n_\Gamma+q+\Lambda$)}\notag\\
			&\leq C_2R^{C_3}+C\det(\Lambda)\qquad\quad \mbox{(by \Cref{lem:lattice})}.
		\end{align}
		
		It remains to estimate the last right hand side term in \eqref{equ-w2}. This can be done as follows.
		
		\medskip

	\textbf{Claim. }Suppose that $\Gamma\in\mathbb{A}_\Lambda$ is a $s$-dimensional sublattice in $\Z^d$. If there exists $\alpha\in \mathscr{I}_\flat(\Gamma;R)$, then
	\begin{equation}
		\label{eq:estimatedeterminant}
		\det(\Lambda)\leq R^{sC_1}.
	\end{equation}
	
	Indeed since $\alpha\in \mathscr{I}_\flat(\Gamma;R)$, there exists $Z_\Gamma^\alpha$ such that $\dim\Gamma^\alpha=\dim\Gamma=s$. Thus, we can choose $s+1$ independent points $\{x^0, x^1, \dots, x^s\} \subset Z_\Gamma^\alpha$. By defining the difference vectors $\mathbf{w}^j = x^j-x^0$ for $j=1,\dots,s$, we obtain $s$ linearly independent lattice vectors in $Z_\Gamma^\alpha$. Obviously,
	\[
	|\mathbf{w}^j|\leq \text{diam}(Z_\Gamma^\alpha),\quad\,j=1,\cdots,s.
	\]
	Let $\tilde{\Lambda}=\text{span}_\Z\{\mathbf{w}^1,\cdots,\mathbf{w}^s\}$. Then $\tilde{\Lambda}$ is a sublattice of $\Lambda$, which implies (see \cite[section I.2.2]{Cas97}):
	\begin{equation}\label{equ1}
		\det(\Lambda)\leq \det(\tilde{\Lambda}).
	\end{equation}
	By Hadamard's inequality and \eqref{equ-Cluster-diam}, we deduce that
	\begin{equation}\label{equ2}
		\det(\tilde{\Lambda})\leq \prod_{j=1}^s|\mathbf{w}^j|\leq \text{diam}(Z_\Gamma^\alpha)^s\leq \text{diam}(Q_\Gamma^\alpha)^s \leq R^{sC_1}.
	\end{equation}
	By $\eqref{equ1}$ and \eqref{equ2}, we obtain \eqref{eq:estimatedeterminant}.

		\medskip
		
		It then follows from \eqref{equ-w2} and \eqref{eq:estimatedeterminant} that \eqref{eq:estimationsupport} holds.
\end{proof}
	
	\section{Quantitative observability for the free Schr\"{o}dinger equation}
	
	The goal of this part is to prove \Cref{Thm1}.

	\subsection{Weak observability}
	In this subsection, we aim to establish the weak observability estimate (see \Cref{prop-weak-obs} below). Analogous to Proposition 3.6 in \cite{BZ25}, we provide a quantitative weak observability inequality for initial data $u_0\in L^2_{n_\Gamma+\Gamma}$, where $n_\Gamma\in \Lambda^\perp$ is given by \Cref{lem:lattice}.
	
	\begin{pr}\label{prop-weak-obs}
		Let $0<T<2\pi$, $\omega=\prod_{j=1}^{d}(a_j,b_j)\subset\T^d$ be an open set.
		Assume that the observability \eqref{equ-obs} with $C^{(s)}_{\text{obs}}(T,\omega)$ holds  for all initial data with Fourier supports of dimensions $\leq s\leq d-1$. 
		If $\Gamma\in\mathbb{A}_\Lambda$ and $\dim \Gamma=s+1$,
		then for all $u_0\in L^2_{n_\Gamma+\Gamma}$,
		\begin{equation}\label{equ-weakobs-lowerdim}
			\|u_0\|_{L^2(\T^d)}^2\leq 2C^{(s)}_{\text{obs}}(T,\omega) \int_0^T\|e^{-it\Delta}u_0\|_{L^2(\omega)}^2dt+ 2\|\Pi_{B_\Gamma} u_0\|_{L^2(\T^d)}^2.
		\end{equation}
		In \eqref{equ-weakobs-lowerdim},
		\begin{equation*}
			B_\Gamma\subset B(0,c_1(R_{s+1})^{c_2}),
		\end{equation*}
		where  $c_1,c_2>0$ depend only on $d,s$, and $R_{s+1}=\big(4C_0(2^d+1)C_{\text{obs}}^{(s)}(T,\omega)\big)^4$ with $C_0$ be given by \Cref{lem-psi}.
		
	\end{pr}
	
	The proof of \Cref{prop-weak-obs} is based on mathematical induction, and is reported at the end of this part. It is important to notice that roughly speaking \Cref{prop-weak-obs} states that observability in low dimension implies weak observability for initial data with Fourier support in $n_{\Gamma} + \Gamma$. In particular, starting from this result, one cannot iterate it in a direct way to obtain the (weak) observability for every initial data.
	
	\medskip
	
	We start by comparing in a quantitative way the difference between the observation on $\omega$ with its approximation term involving the band limited function $\psi$.

	\begin{lm}
		\label{lem:errorSchro}
		For every $u_0 \in L^2(\T^d)$, we have 
		\begin{equation}\label{equ-error}
			\int_{\T_t} \int_{\T^d} |\chi_{\omega} \chi_{\omega_t} -\psi||e^{-it\Delta}u_0|^2dxdt\leq (2^d+1)C_0R^{-1/4}\|u_0\|_{L^2(\T^d)}^2,
		\end{equation}
	\end{lm}
	
	\begin{proof}
		We prove \eqref{equ-error} by induction on the dimension $d$ of the torus. 
		
		If $d=1$, by Parseval's equality and \Cref{lem-psi}, we derive that
		\begin{align*}
			&\quad\int_{\T_t}\int_{\T} |\chi_{\omega} \chi_{\omega_t} -\psi||e^{-it\Delta}u_0|^2dxdt\\
			&\leq\int_{\T_t}\int_{\T} |\chi_{\omega}-\varphi_1(x)||e^{-it\Delta}u_0|^2dxdt+ \int_{\T_t}\int_{\T} |\varphi_1(x)\chi_{\omega_t}-\varphi_1(x)\varphi_2(t)||e^{-it\Delta}u_0|^2dxdt\\
			&\leq \int_{\T} |\chi_{\omega}-\varphi_1(x)|\int_{\T_t}\big|\sum_{n\in\Z}c_ne^{i(n\cdot x+|n|^2t)}\big|^2dtdx+\int_{\T_t}|\chi_{\omega_t}-\varphi_2(t)|\int_{\T}|e^{-it\Delta}u_0|^2dxdt\\
			&\leq \sum_{m\in\N}\int_{\T} |\chi_{\omega}-\varphi_1(x)|\big|\sum_{n\in\Z:\,|n|^2=m}c_ne^{in\cdot x}\big|^2dx+\|\chi_{\omega_t}-\varphi_2\|_{L^1(\T)}\|u_0\|_{L^2(\T)}^2\\
			&\leq \sum_{m\in\N}\int_{\T} |\chi_{\omega}-\varphi_1(x)|\big|c_{\sqrt{m}} e^{i {\sqrt{m}}x} + c_{-\sqrt{m}} e^{-i{\sqrt{m}}x} |^2 dx+\|\chi_{\omega_t}-\varphi_2\|_{L^1(\T)}\|u_0\|_{L^2(\T)}^2\\
			&\leq 2\|\chi_{\omega}-\varphi_1\|_{L^1(\T)}\sum_{m\in\N}\sum_{n\in\Z:\,|n|^2=m}|c_n|^2+\|\chi_{\omega_t}-\varphi_2\|_{L^1(\T)}\|u_0\|_{L^2(\T)}^2\\
			&\leq \big(2\|\chi_{\omega}-\varphi_1\|_{L^1(\T)}+\|\chi_{\omega_t}-\varphi_2\|_{L^1(\T)}\big) \|u_0\|_{L^2(\T)}^2\\
			&\leq 3C_0R^{-1/4}\|u_0\|_{L^2(\T)}^2.
		\end{align*}
		Thus, \eqref{equ-error} holds for $d=1$. 
		
		Next, we assume that \eqref{equ-error} holds for $d=k-1$  with $k\geq 2$. For $d=k$, let $\tilde{\psi}(x',t)=(\prod_{j=1}^{k-1}\varphi_j(x_j))\varphi_{k+1}(t)$ with $x'=(x_1,\cdots,x_{k-1})$, and $\omega'=\prod_{j=1}^{k-1}(a_j,b_j),\,\omega_k=(a_k,b_k)$, we then have by the induction property and the previous computation for $d=1$,
		\begin{align*}
			&\quad\int_{\T_t}\int_{\T^k} |\chi_{\omega} \chi_{\omega_t}-\psi||e^{-it\Delta}u_0|^2dxdt\\
			&\leq \int_{\T_t}\int_{\T^k} |\chi_{\omega}\chi_{\omega_t}-\tilde{\psi}\chi_{\omega_k}||e^{-it\Delta}u_0|^2dxdt+ \int_{\T_t}\int_{\T^k} |\tilde{\psi}\chi_{\omega_k}-\psi||e^{-it\Delta}u_0|^2dxdt\\
			&\leq \int_{\T_t}\int_{\T^{k-1}}|\chi_{\omega'}\chi_{\omega_t}-\tilde{\psi}|\Big(\int_{\T_{x_k}}\Big|\sum_{n_k\in\Z}\big(\sum_{n'\in\Z^{k-1}}c_{n',n_k}e^{i(n'\cdot x'+|n'|^2t)}\big)e^{in_kx_k+in_k^2t}\Big|^2dx_k\Big)dx'dt\\
			&\quad+ \int_{\T_t}\int_{\T_{x_k}}|\chi_{\omega_k}-\varphi_k|\Big(\int_{\T^{k-1}}\Big|\sum_{n'}\big(\sum_{n_k}c_{n',n_k}e^{i(n_kx_k+n_k^2t)}\big)e^{i(n'\cdot x'+|n'|^2t)}\Big|^2dx'\Big)dx_kdt\\
			&\leq \sum_{n_k\in\Z} \int_{\T_t}\int_{\T^{k-1}}|\chi_{\omega'}\chi_{\omega_t}-\tilde{\psi}|\big|\sum_{n'\in\Z^{k-1}}c_{n',n_k}e^{i(n'\cdot x'+|n'|^2t)}\big|^2dx'dt\\
			&\quad+\sum_{n'\in\Z^{k-1}}\int_{\T_t}\int_{\T_{x_k}}|\chi_{\omega_k}-\varphi_k|\big|\sum_{n_k\in\Z}c_{n',n_k}e^{i(n_kx_k+n_k^2t)}\big|^2dx_kdt\\
			&\leq \sum_{n_k\in\Z}\Big((2^{k-1}+1)C_0R^{-1/4}\sum_{n'\in\Z^{k-1}}|c_{n',n_k}|^2\Big)+2\|\chi_{\omega_k}-\varphi_k\|_{L^1(\T)}\|u_0\|_{L^2(\T^k)}^2\\
			&\leq (2^k+1)C_0R^{-1/4}\|u_0\|_{L^2(\T^k)}^2.
		\end{align*}
		Therefore, \eqref{equ-error} holds for $d=k$ and we obtain \eqref{equ-error} by induction.	
	\end{proof}
	
	\begin{rmk}
		A property that is hidden and used in the above proof is that for the Schrödinger equation \eqref{eq:schro} in one-dimensional torus $\T$, the following Strichartz estimate holds
		\begin{equation}
			\label{eq:Strichartz1D}
			\|e^{-it \Delta} u_0\|_{L_x^{\infty}(\T;L^2(0,T))} \leq C \|u_0\|_{L^2(\T)}.
		\end{equation}
	\end{rmk}

	\begin{proof}[Proof of \Cref{prop-weak-obs}]
		Assume that $\Gamma=q+\Lambda\in\mathbb A_\Lambda$
		with $\dim\Gamma=s+1$, and set $\widetilde\Gamma:=n_\Gamma+\Gamma$.
		Let $u_0\in L^2_{\widetilde\Gamma}$.
		For simplicity, we write $\widetilde\Gamma$ as $\Gamma$ in the sequel.
		We define the cluster components of $u_0$ by
		\[
		u_{0,\alpha}(x) := \sum_{n\in  Z_\Gamma^\alpha} c_{n} e^{in\cdot x}, 
		\quad
		u_{0}(x) = \sum_{\alpha\in\mathscr{I}(\Gamma;R)} u_{0,\alpha}(x),
		\]
		and let $u_\alpha=e^{-it\Delta}u_{0,\alpha}$ denote the corresponding solutions of \eqref{eq:schro}. 
		To prove \eqref{equ-weakobs-lowerdim} for $\mathrm{dim}\ \Gamma=s+1$, we proceed as follows.
		
		\medskip
		
		Using the cluster decomposition and the inductive hypothesis, it follows that
		\begin{align}\label{equ-boundwithR}
			\quad\|u_0\|_{L^2(\T^d)}^2&=\sum_{\alpha\in \mathscr{I}(\Gamma;R)}\sum_{n\in  Z_\Gamma^\alpha}|c_n|^2\notag\\
			&=\sum_{\alpha\in \mathscr{I}_\sharp(\Gamma;R)}\|u_{0,\alpha}\|_{L^2(\T^d)}^2+\sum_{\alpha\in \mathscr{I}_\flat(\Gamma;R)}\|u_{0,\alpha}\|_{L^2(\T^d)}^2\notag\\
			&\leq C_{\text{obs}}^{(s)}(T,\omega) \sum_{\alpha\in \mathscr{I}_\sharp(\Gamma;R)}\int_{0}^{T} \|e^{-it\Delta}u_{0,\alpha}\|_{L^2(\omega)}^2dt+\sum_{\alpha\in \mathscr{I}_\flat(\Gamma;R)}\|u_{0,\alpha}\|_{L^2(\T^d)}^2\notag\\
			&\hspace{6cm} \mbox{(since $\mathrm{dim}\ \Gamma^\alpha\leq s$ for $\alpha\in \mathscr{I}_\sharp(\Gamma;R)$)}\notag\\
			&\leq  C_{\text{obs}}^{(s)}(T,\omega) \sum_{\alpha\in \mathscr{I}(\Gamma;R)}\int_{0}^{T} \|e^{-it\Delta}u_{0,\alpha}\|_{L^2(\omega)}^2dt+\sum_{\alpha\in \mathscr{I}_\flat(\Gamma;R)}\|u_{0,\alpha}\|_{L^2(\T^d)}^2\notag\\
			&=C_{\text{obs}}^{(s)}(T,\omega)\big( \int_{0}^{T}\|e^{-it\Delta}u_0\|_{L^2(\omega)}^2dt +\mathcal{R} \big)+\sum_{\alpha\in \mathscr{I}_\flat(\Gamma;R)}\|u_{0,\alpha}\|_{L^2(\T^d)}^2,
		\end{align}
		where 
		\begin{equation*}
			\mathcal{R}=\sum_{\alpha\in\mathscr{I}(\Gamma;R)} \int_{0}^{T} \|e^{-it\Delta}u_{0,\alpha}\|_{L^2(\omega)}^2dt- \int_{0}^{T}\|e^{-it\Delta}u_0\|_{L^2(\omega)}^2dt.
		\end{equation*}
		We then estimate $\mathcal{R}$. Let $\psi$ be given in \eqref{eq:defpsi}. By the support property of $\psi$ recalled in \eqref{eq:supppsi} and the separation between clusters recalled in \eqref{equ-Cluster-dist}, we get the following orthogonality property
		\[
		\sum_{\alpha\in\mathscr{I}(\Gamma;R)} \int_{\T_t} \int_{\T^d} | u_{\alpha} |^2 \psi dx dt=\int_{\T_t} \int_{\T^d} |\sum_{\alpha\in\mathscr{I}(\Gamma;R)} u_{\alpha} |^2 \psi dx dt.
		\]
		Thus, by the error estimate \eqref{equ-error} of \Cref{lem:errorSchro}, we derive that
		\begin{align}\label{equ-upperR}
			|\mathcal{R}|&\leq \sum_{\alpha\in\mathscr{I}(\Gamma;R)}\int_{\T_t} \int_{\T^d}|\chi_{\omega} \chi_{\omega_t}-\psi||u_\alpha|^2 dxdt+ \int_{\T_t}\int_{\T^d} |\chi_{\omega} \chi_{\omega_t}-\psi||u|^2dxdt\notag\\
			&\leq (2^d+1)C_0R^{-1/4}\big(\sum_{\alpha\in\mathscr{I}(\Gamma;R)}\|u_{0,\alpha}\|_{L^2(\T^d)}^2+\|u_0\|_{L^2(\T^d)}^2\big)\notag\\
			&=2(2^d+1)C_0R^{-1/4}\|u_0\|_{L^2(\T^d)}^2.
		\end{align}
		Since $R_{s+1}=\big(4C_0(2^d+1)C_{\text{obs}}^{(s)}(T,\omega)\big)^4$, taking $R= R_{s+1} $, we get
		\[
		C_{\text{obs}}^{(s)}(T,\omega)2(2^d+1)C_0R^{-1/4}=\frac{1}{2},
		\]
		which, combining with \eqref{equ-boundwithR}, \eqref{equ-upperR}, implies
		\begin{equation}\label{equ-w3}
			\|u_0\|_{L^2(\T^d)}^2\leq 2C_{\text{obs}}^{(s)}(T,\omega) \int_{0}^{T}\|e^{-it\Delta}u_0\|_{L^2(\omega)}^2dt+2\sum_{\alpha\in \mathscr{I}_\flat(\Gamma;R)}\|u_{0,\alpha}\|_{L^2(\T^d)}^2.
		\end{equation}

		Finally, to estimate $\sum_{\alpha\in \mathscr{I}_\flat(\Gamma;R)}\|u_{0,\alpha}\|_{L^2(\T^d)}^2$, we use the estimate on the support \eqref{eq:estimationsupport} of \Cref{lem:support}. Hence for $u_{0,\flat}:=\sum_{\alpha\in\mathscr{I}_\flat(\Gamma,R)}u_{0,\alpha}$,
		\[
		\text{supp }\widehat{u_{0,\flat}}\subset B(0,c_1R^{c_2}).
		\]
		This, together with \eqref{equ-w3}, implies \eqref{equ-weakobs-lowerdim}.
		Therefore, the proof is complete. 
	\end{proof}

	\subsection{Low-frequency part}
	In this section, we build the quantitative observability for the low-frequency part.
	
	\subsubsection{Observability for the heat equation}
	We first give the quantitative observability for the $d$-dimensional homogeneous heat equation:
	\begin{equation}\label{equ-heat}
		\partial_t u-\Delta u=0,\qquad u(x,0)=u_0\in L^2(\T^d).
	\end{equation} 
	\begin{pr}\label{prop:Heat}
		Let $\omega\subset\T^d$ be a measurable set of positive Lebesgue measure $|\omega|>0$. Then there exists a positive constant $C>0$ such that for any $T>0$ and any solution $u$ of \eqref{equ-heat}, 
		\begin{equation}
			\label{eq:obsheat}
			\|u(T)\|_{L^2(\T^d)}^2\leq  C \exp\left(C \frac{(\log(|\omega|))^2}{T}\right)\int_{0}^T\|u(t)\|_{L^2(\omega)}^2dt.
		\end{equation}
	\end{pr} 
	It is interesting to note that the estimate \eqref{eq:obsheat} matches with the one given in \cite[Theorem 1.3]{LL21} for $\omega = B(x_0,r)$.
	
	\begin{proof}
		The proof is based on two steps. We first establish a spectral estimate quantitative in function of the measure of $\omega$, based on Tur\'{a}n-Nazarov's inequality. From this estimate, we employ the by now classical time iteration Lebeau-Robbiano strategy to deduce the associated observability estimate for the heat equation.
		
		\medskip
		
		\emph{Step 1:  Spectral estimate.} Let $\omega\subset\T^d$ with $|\omega|>0$, then for any $N>0$ and every $f\in L^2(\T^d)$,
		\begin{equation}\label{equ-spec}
			\|\Pi_N f\|_{L^2(\T^d)}^2\leq \left(\frac{28d}{|\omega|}\right)^{4(N+1)d} \int_{\omega}|\Pi_N f|^2dx.
		\end{equation}
		
		Indeed, take $f=\Pi_N f$. By the high-dimensional version of Tur\'{a}n-Nazarov's inequality (see \cite[Theorem 2]{FM06}), we get
		\begin{equation}\label{equ-SmallnessSup}
			\sup_{\T^d}|f|\leq \left(\frac{14d}{|\omega|}\right)^{(2N+1)d}\sup_{\omega}|f|.
		\end{equation}
		
		Next, we follow the proof of the Corollary in \cite{Kov01}.
		From \eqref{equ-SmallnessSup}, we have 
		\begin{equation*}
			\big|\big\{x\in \T^d:\,|f(x)|< \big(\frac{\eps}{14d}\big)^{(2N+1)d}\|f\|_{L^\infty(\T^d)}\big\}\big|\leq \eps
		\end{equation*}
		for any $\eps>0$. Let $\eps=\frac{|\omega|}{2}$, we get
		\begin{equation*}
			\big|\big\{x\in \T^d:\,|f(x)|< \big(\frac{|\omega|}{28d}\big)^{(2N+1)d}\|f\|_{L^\infty(\T^d)}\big\}\big|\leq \frac{|\omega|}{2}.
		\end{equation*}
		Thus,
		\begin{align*}
			\int_\omega|f|^2dx&\geq \int_\omega\chi_{|f|\ge (\frac{|\omega|}{28d})^{(2N+1)d} \|f\|_{L^\infty(\T^d)}}\cdot|f|^2\\
			&\geq \big(\frac{|\omega|}{28d}\big)^{2(2N+1)d}\cdot \|f\|_{L^\infty(\T^d)}^2\cdot \frac{|\omega|}{2}\\
			&\geq \big(\frac{|\omega|}{28d}\big)^{4(N+1)d}\cdot \|f\|_{L^2(\T^d)}^2,
		\end{align*}
		which concludes the proof of \eqref{equ-spec}.

		\medskip
		
		\emph{Step 2: Quantitative Lebeau-Robbiano strategy.} We now employ Lebeau-Robbiano's strategy, see \cite{LR95}, starting from the dissipation of the heat-semigroup
		\begin{equation}
			\label{eq:dissipationheat}
			\|e^{t\Delta}(1-\Pi_N)f\|_{L^2(\T^d)}^2=\sum_{n\in\Z^d:|n|>N}|c_n|^2 e^{-2t|n|^2}
			\leq e^{-2 t N^2}\|f\|_{L^2(\T^d)}^2,\quad 0  < t < \infty.
		\end{equation}
		We can use \cite[Lemma 2.6]{LL21} directly, which is borrowed from \cite{Mil10}. The combination of \eqref{equ-spec} and \eqref{eq:dissipationheat}	leads to \eqref{eq:obsheat}.
	\end{proof}
	
	Next, we consider the $d$-dimensional inhomogeneous backward heat equation
	\begin{equation}\label{equ-InhomoHeatBack}
		\partial_t u+(\Delta-1) u=F(t,x)\in L^2((0,\infty); L^2(\T^d)),\quad u(x,T)=u_T\in H^2(\T^d).
	\end{equation} 
	\begin{cor}\label{cor:ObsHeatBack}
		There exists a constant $C>1$, depending only on $d$, such that for any $T>0$ and any solution  $u$ of \eqref{equ-InhomoHeatBack} we have
		\[
		\|u(0)\|_{L^2(\T^d)}^2\leq C \exp\left(C \frac{(\log(|\omega|))^2}{T}\right)\int_0^T\big(\|(-\Delta+1)u(t)\|_{L^2(\omega)}^2+\|F(t)\|_{L^2(\T^d)}^2\big)dt.
		\]
	\end{cor}
	\begin{proof}
		By reversing the time, i.e. set $t$ to $T-t$, it is equivalent to consider
		\begin{equation*}
			\partial_t u+(-\Delta+1) u=F(t,x)\in L^2((0,\infty); L^2(\T^d)),\quad u(x,0)=u_0\in H^2(\T^d).
		\end{equation*}
		and obtain
		\begin{align*}
			\|u(T)\|_{L^2(\T^d)}^2\notag
			\leq C \exp\left(C \frac{(\log(|\omega|))^2}{T}\right)\int_0^T\big(\|(-\Delta+1)u(t)\|_{L^2(\omega)}^2+\|F(t)\|_{L^2(\T^d)}^2\big)dt.
		\end{align*}
		We then decompose the solution into two parts
		\begin{equation*}
			\partial_t u_1+(-\Delta+1) u_1=0 ,\quad u_1(x,0)=u_0\in H^2(\T^d),
		\end{equation*}
		and 
		\begin{equation*}
			\partial_t u_2+(-\Delta+1) u_2=F ,\quad u_2(x,0)=0.
		\end{equation*}
		We have by the observability esimate \eqref{eq:obsheat}
		$$ \|u_1(T)\|_{L^2(\T^d)}^2 \leq \|(-\Delta + 1 )u_1(T)\|_{L^2(\T^d) }^2\leq C \exp\left(C \frac{(\log(|\omega|))^2}{T}\right)\int_0^T\|(-\Delta+1)u_1(t)\|_{L^2(\omega)}^2 dt.$$
		We also have by $L^2$ parabolic regularity, see for instance \cite{Eva10}, that
		$$ \|u_2 (T)\|_{L^2(\T^d)}^2 + \int_0^T \|(-\Delta + 1 )u_2(t)\|_{L^2(\T^d) }^2  dt \leq \exp\left(C \frac{(\log(|\omega|))^2}{T}\right)\int_0^T\|F(t)\|_{L^2(\T^d)}^2dt.$$
		The combination of the two inequalities easily gives the result.
	\end{proof}
	
	\subsubsection{Quantitative unique continuation}
	Analogous to \cite[Proposition 5.1]{SSY25}, we give the following quantitative unique continuation
	estimate for the $d$-dimensional Schr\"{o}dinger equation.
	
	\begin{pr}\label{prop:lowfreq}
		There exist some constants $h_0\in(0,1)$ and $C>1$, depending only on $d$, such that for any $T>0,\,0<h<h_0(1+|\omega|^{-1}T^{-3})^{-1}$ and $f\in H^2(\T^d)$, we have
		\begin{align}\label{equ-Quanti-UC}
			\|f\|_{L^2(\T^d)}^2&\leq Ch(1+\frac{1}{T^2})\|(-\Delta+1)f\|_{L^2(\T^d)}^2\notag\\
			&+C\exp\left(\frac{2T^2}{h}+C \frac{(\log(|\omega|))^2}{T}\right) \int_0^T\|(-\Delta+1)e^{it(-\Delta+1)}f\|_{L^2(\omega)}^2dt.
		\end{align}
	\end{pr}
	\begin{proof}
		Following the proof of Proposition 5.1 in \cite{SSY25}, using notably the FBI transform we conclude the following estimate based on \Cref{cor:ObsHeatBack},
		\begin{align*}
			\|f\|_{L^2(\T^d)}^2&\leq Ch (1+\frac{1}{T^2})\|(-\Delta+1)f\|_{L^2(\T^d)}^2\\
			&+\frac{C T^2}{h}\exp\left(\frac{T^2}{h}+C \frac{(\log(|\omega|))^2}{T}\right)\int_{0}^{T}\|(-\Delta+1)e^{it(-\Delta+1)}f(t)\|_{L^2(\omega)}^2dt\\
			&+\left(\frac{C T}{ h}\exp\left(-\frac{3T^2}{h}+C\frac{(\log(|\omega|))^2}{T}\right)+2e^{-\frac{18T^2}{h}}\right)\|f\|_{L^2(\T^d)}^2.
		\end{align*}
		Then, we can take $h_0$ small enough, independent of $T$ and $\omega$, such that for any $0<h<h_0(1+|\omega|^{-1}T^{-3})^{-1}$, we have
		\[
		\frac{C T}{ h}\exp\left(-\frac{3T^2}{h}+C\frac{(\log(|\omega|))^2}{T}\right)+2e^{-\frac{18T^2}{h}}<\frac{1}{2}.
		\]
		This completes the proof of \Cref{prop:lowfreq}.
	\end{proof}
	
	\subsection{Proof of \Cref{Thm1}}
	
	In this part, we aim to
	prove the following theorem which clearly implies \Cref{Thm1}.
	
	\begin{tm}\label{thm:obslowdim}
		If $\Lambda$ is a primitive sublattice with $\dim\Lambda=s\leq d$, $\Gamma\in\mathbb{A}_\Lambda$ and $u_0\in L^2_\Gamma$, then
		\begin{equation}\label{equ-obs-s}
			\|u_0\|_{L^2(\T^d)}^2\leq C^{(s)}_{\text{obs}}(T,\omega)\int_{0}^T\|e^{-it\Delta}u_0\|_{L^2(\omega)}^2dt,
		\end{equation}
		where
		\begin{equation}\label{equ-obs-c0}
			C^{(0)}_{\text{obs}}(T,\omega)=\frac{1}{T|\omega|},
		\end{equation}
		and 
		\begin{equation}\label{equ-obs-cs}
			C^{(s)}_{\text{obs}}(T,\omega)=C_s''\underbrace{\exp\{\exp\{\cdots\exp\{}_{s\text{ times }}\big(\frac{C_s''}{(T|\omega|)^{C_s'}}+C_s''\frac{(\log(|\omega|))^2}{T}\big)\}\}\}\qquad\mbox{ if } 1\leq s\leq d.
		\end{equation}
		In \eqref{equ-obs-cs}, the constants $C_s',C_s''>0$ depend only on $d,s$.
	\end{tm}
	
	\begin{proof}
		If $s = 0$, we have that $u_0 = c_n e^{i n \cdot x}$ for some $n \in \Z^d$ so $e^{-i t \Delta} u_0 = c_n e^{i n \cdot x} e^{i|n|^2 t}$, then 
		\begin{equation}\label{equ-s0}
			\int_{0}^{T} \|e^{-it\Delta}u_0\|_{L^2(\omega)}^2dt = T |\omega| |c_n|^2 = T |\omega| \|u_0\|_{L^2(\T^d)}^2,
		\end{equation}
		so \eqref{equ-obs-s} holds for $s=0$. Remark that this basis induction case is similar to \cite[Lemma 3.5]{BZ25}.
		
		If $s=1$, let $\tilde{u_0}=e^{in_\Gamma\cdot x}u_0\in L^2_{n_\Gamma+\Gamma}$ with $n_\Gamma\in\Lambda^\perp$ be given by \Cref{lem:lattice}, it follows from \Cref{prop-weak-obs} that
		\begin{equation}\label{equ-s1}
			\|\tilde{u_0}\|_{L^2(\T^d)}^2\leq 2C^{(0)}_{\text{obs}}(T,\omega) \int_0^T\|e^{-it\Delta}\tilde{u_0}\|_{L^2(\omega)}^2dt+ 2\|\Pi_{B_\Gamma} \tilde{u_0}\|_{L^2(\T^d)}^2,
		\end{equation}
		where 
		\begin{equation}\label{equ-s2}
			B_\Gamma\subset B(0,c_1(R_{1})^{c_2}),\qquad R_{1}=\big(4C_0(2^d+1)C_{\text{obs}}^{(0)}(T,\omega)\big)^4.
		\end{equation}
		We denote $\tilde{U_0}=(-\Delta+1)^{-1}\tilde{u_0}\in H^2(\T^d)$ and $\tilde{u}(t)=e^{it(-\Delta+1)}\tilde{u_0}$.
		By the quantitative unique continuation \eqref{equ-Quanti-UC}, we have
		\begin{equation}\label{equ-U0}
			\|\tilde{U_0}\|_{L^2(\T^d)}^2\leq Ch(1+\frac{1}{T^2})\|\tilde{u_0}\|_{L^2(\T^d)}^2+C\exp\left(\frac{2T^2}{h}+C \frac{(\log(|\omega|))^2}{T}\right)\int_{0}^{T}\|\tilde{u}(t)\|_{L^2(\omega)}^2dt,
		\end{equation}
		for any $T>0,\,0<h<h_0(1+|\omega|^{-1}T^{-3})^{-1}$.
		From the definition of $\tilde{U_0}$ and \eqref{equ-s2}, we obtain
		\begin{align*}
			\|\Pi_{B_\Gamma} \tilde{u_0}\|_{L^2(\T^d)}^2&\leq\sum_{n\in\Z^d:\,|n|\leq c_1(R_{1})^{c_2}}|c_n|^2=\sum_{|n|\leq c_1(R_{1})^{c_2}}(|n|^2+1)^2(|n|^2+1)^{-2}|c_n|^2\\
			&\leq (c_1^2(R_{1})^{2c_2}+1)^2\sum_{|n|\leq c_1(R_{1})^{c_2}}(|n|^2+1)^{-2}|c_n|^2\leq 4c_1^4(R_{1})^{4c_2}\|\tilde{U_0}\|_{L^2(\T^d)}^2.
		\end{align*}
		This, combining with \eqref{equ-s1}, \eqref{equ-U0}, deduces that
		\begin{align*}
			\|\tilde{u_0}\|_{L^2(\T^d)}^2&\leq 2C^{(0)}_{\text{obs}}(T,\omega) \int_0^T\|e^{-it\Delta}\tilde{u_0}\|_{L^2(\omega)}^2dt+8c_1^4(R_{1})^{4c_2}\|\tilde{U_0}\|_{L^2(\T^d)}^2\\
			&\leq 2C^{(0)}_{\text{obs}}(T,\omega) \int_0^T\|e^{-it\Delta}\tilde{u_0}\|_{L^2(\omega)}^2dt+ C(R_{1})^{4c_2}h(1+\frac{1}{T^2})\|\tilde{u_0}\|_{L^2(\T^d)}^2\\
			&\hspace{1cm}+C(R_{1})^{4c_2}\exp\left(\frac{2T^2}{h}+C\frac{(\log(|\omega|))^2}{T}\right)\int_{0}^{T}\|\tilde{u}(t)\|_{L^2(\omega)}^2dt.
		\end{align*} 
		Since $0<|\omega|<1$ and $0<T<1$,
		taking $h=\varepsilon(1+(|\omega| T)^{-C'})^{-1}$ with $\varepsilon$ small enough to absorb the second right hand side in the previous inequality, by \eqref{equ-obs-c0} and \eqref{equ-s2}, we obtain
		\begin{equation*}
			\|\tilde{u_0}\|_{L^2(\T^d)}^2\leq   C\exp\big(\frac{C}{(T|\omega|)^{C'}}+C\frac{(\log(|\omega|))^2}{T}\big)\int_0^T\|e^{-it\Delta}\tilde{u_0}\|_{L^2(\omega)}^2dt.
		\end{equation*}
		Since $n_\Gamma\in\Lambda^\perp$,
		by applying \Cref{prop:equiv}, we conclude \eqref{equ-obs-s} for $s=1$.
		
		Next, we assume that \eqref{equ-obs-s} holds for $s\leq\nu-1$ with $\nu\ge2$. For $s=\nu$, we proceed similarly to the case $s=1$ and obtain
		\begin{align}\label{equ-s3}
			\|\tilde{u_0}\|_{L^2(\T^d)}^2
			&\leq 2C^{(\nu-1)}_{\text{obs}}(T,\omega) \int_0^T\|e^{-it\Delta}\tilde{u_0}\|_{L^2(\omega)}^2dt+ C(R_{\nu})^{4c_2}h(1+\frac{1}{T^2})\|\tilde{u_0}\|_{L^2(\T^d)}^2\notag\\
			&\hspace{1cm}+C(R_{\nu})^{4c_2}\exp\left(\frac{2T^2}{h}+C\frac{(\log(|\omega|))^2}{T}\right)\int_{0}^{T}\|\tilde{u}(t)\|_{L^2(\omega)}^2dt.
		\end{align} 
		Recall that $R_\nu=\big(4C_0(2^d+1)C_{\text{obs}}^{(\nu-1)}(T,\omega)\big)^4$, where $C_{\text{obs}}^{(\nu-1)}(T,\omega)$ satisfies \eqref{equ-obs-cs},
		taking 
		\[
		h=\eps\Big(\underbrace{\exp\{\exp\{\cdots\exp\{}_{\nu-1\text{ times }}\big(\frac{C''}{(T|\omega|)^{C'}}+C''\frac{(\log(|\omega|))^2}{T}\big)\}\}\}\Big)^{-1}
		\]
		with $\varepsilon$ small enough, it follows from \eqref{equ-s3} and \Cref{prop:equiv} that \eqref{equ-obs-s} holds for $s=\nu$. Thus, the proof of \Cref{thm:obslowdim} is completed by induction.
	\end{proof}

	\bibliographystyle{alpha}
	\bibliography{SchrodingerTd}

\begin{thebibliography}{BGSZ26}

\bibitem[AL25]{AL25}
Paul Alphonse and David Lafontaine.
\newblock Observability estimates for the schr\"odinger equation on the
  equilateral triangle.
\newblock {\em arXiv eprint: 2509.24642}, 2025.

\bibitem[AM14]{AM14}
Nalini Anantharaman and Fabricio Maci\`a.
\newblock Semiclassical measures for the {S}chr\"{o}dinger equation on the
  torus.
\newblock {\em J. Eur. Math. Soc. (JEMS)}, 16(6):1253--1288, 2014.

\bibitem[AT25]{AT25}
Paul Alphonse and Nikolay Tzvetkov.
\newblock A smoothing effect for the fractional {S}chr\"odinger equations on
  the circle and observability.
\newblock {\em Int. Math. Res. Not. IMRN}, (8):Paper No. rnaf100, 20, 2025.

\bibitem[BBZ13]{BBZ13}
Jean Bourgain, Nicolas Burq, and Maciej Zworski.
\newblock Control for {S}chr\"{o}dinger operators on 2-tori: rough potentials.
\newblock {\em J. Eur. Math. Soc. (JEMS)}, 15(5):1597--1628, 2013.

\bibitem[BGSZ26]{BGSZ26}
Nicolas Burq, Pierre Germain, Massimo Sorella, and Hui Zhu.
\newblock Trace and observability inequalities for {Laplace} eigenfunctions on
  the torus.
\newblock {\em Forum of Mathematics, Sigma}, 14:e47, 2026.

\bibitem[Bou14]{Bou14}
Jean Bourgain.
\newblock On the control problem for {S}chr\"{o}dinger operators on tori.
\newblock In {\em Geometric aspects of functional analysis}, volume 2116 of
  {\em Lecture Notes in Math.}, pages 97--105. Springer, Cham, 2014.

\bibitem[BZ04]{BZ04}
Nicolas Burq and Maciej Zworski.
\newblock Geometric control in the presence of a black box.
\newblock {\em J. Amer. Math. Soc.}, 17(2):443--471, 2004.

\bibitem[BZ12]{BZ12}
Nicolas Burq and Maciej Zworski.
\newblock Control for {S}chr\"{o}dinger operators on tori.
\newblock {\em Math. Res. Lett.}, 19(2):309--324, 2012.

\bibitem[BZ19]{BZ19}
Nicolas Burq and Maciej Zworski.
\newblock Rough controls for {S}chr\"{o}dinger operators on 2-tori.
\newblock {\em Ann. H. Lebesgue}, 2:331--347, 2019.

\bibitem[BZ25]{BZ25}
Nicolas Burq and Hui Zhu.
\newblock Observability of schr{\"o}dinger propagators on tori in rough
  settings.
\newblock {\em arXiv e-prints: 2509.23965}, 2025.

\bibitem[Cas97]{Cas97}
J.~W.~S. Cassels.
\newblock {\em An introduction to the geometry of numbers}.
\newblock Classics in Mathematics. Springer-Verlag, Berlin, 1997.
\newblock Corrected reprint of the 1971 edition.

\bibitem[Cor07]{Cor07}
Jean-Michel Coron.
\newblock {\em Control and nonlinearity}, volume 136 of {\em Mathematical
  Surveys and Monographs}.
\newblock American Mathematical Society, Providence, RI, 2007.

\bibitem[DJN22]{DJN22}
Semyon Dyatlov, Long Jin, and St{\'e}phane Nonnenmacher.
\newblock Control of eigenfunctions on surfaces of variable curvature.
\newblock {\em Journal of the American Mathematical Society}, 35(2):361--465,
  2022.

\bibitem[DSS24]{DSS24}
Art\=uras Dubickas, Min Sha, and Igor~E. Shparlinski.
\newblock Euclidean minima of algebraic number fields.
\newblock {\em Arch. Math. (Basel)}, 122(4):405--414, 2024.

\bibitem[Eva10]{Eva10}
Lawrence~C. Evans.
\newblock {\em Partial Differential Equations}, volume~19 of {\em Graduate
  Studies in Mathematics}.
\newblock American Mathematical Society, Providence, RI, 2 edition, 2010.

\bibitem[FM06]{FM06}
Natacha Fontes-Merz.
\newblock A multidimensional version of turán's lemma.
\newblock {\em Journal of Approximation Theory}, 140(1):27--30, 2006.

\bibitem[Jaf90]{Jaf90}
Stéphane Jaffard.
\newblock Contr\^{o}le interne exact des vibrations d'une plaque rectangulaire.
\newblock {\em Portugal. Math.}, 47(4):423--429, 1990.

\bibitem[Kah62]{Kah62}
Jean-Pierre Kahane.
\newblock Pseudo-périodicité et séries de fourier lacunaires.
\newblock {\em Annales scientifiques de l'École Normale Supérieure},
  79(2):93--150, 1962.

\bibitem[KL05]{KL05}
Vilmos Komornik and Paola Loreti.
\newblock {\em Fourier series in control theory}.
\newblock Springer Monographs in Mathematics. Springer-Verlag, New York, 2005.

\bibitem[Kov01]{Kov01}
Oleg Kovrijkine.
\newblock Some results related to the logvinenko-sereda theorem.
\newblock {\em Proceedings of the American Mathematical Society},
  129(10):3037--3047, 2001.

\bibitem[LBM23]{LBM23}
K\'evin Le~Balc'h and J\'er\'emy Martin.
\newblock Observability estimates for the schr{\"o}dinger equation in the plane
  with periodic bounded potentials from measurable sets.
\newblock {\em arXiv eprint: 2304.08050}, 2023.

\bibitem[Lis14]{Lis14}
Pierre Lissy.
\newblock On the cost of fast controls for some families of dispersive or
  parabolic equations in one space dimension.
\newblock {\em SIAM Journal on Control and Optimization}, 52(4):2651--2676,
  2014.

\bibitem[LL21]{LL21}
Camille Laurent and Matthieu L\'eautaud.
\newblock Observability of the heat equation, geometric constants in control
  theory, and a conjecture of {L}uc {M}iller.
\newblock {\em Anal. PDE}, 14(2):355--423, 2021.

\bibitem[LR95]{LR95}
Gilles Lebeau and Luc Robbiano.
\newblock Contr\^ole exact de l'\'equation de la chaleur.
\newblock {\em Comm. Partial Differential Equations}, 20(1-2):335--356, 1995.

\bibitem[McM84]{McM84}
Peter McMullen.
\newblock Determinants of lattices induced by rational subspaces.
\newblock {\em Bull. London Math. Soc.}, 16(3):275--277, 1984.

\bibitem[Mil04]{Mil04}
Luc Miller.
\newblock How violent are fast controls for schr{\"o}dinger and plate
  vibrations?
\newblock {\em Archive for Rational Mechanics and Analysis}, 172(3):429--456,
  2004.

\bibitem[Mil05]{Mil05}
Luc Miller.
\newblock Controllability cost of conservative systems: resolvent condition and
  transmutation.
\newblock {\em Journal of Functional Analysis}, 218(2):425--444, 2005.

\bibitem[Mil10]{Mil10}
Luc Miller.
\newblock A direct {L}ebeau-{R}obbiano strategy for the observability of
  heat-like semigroups.
\newblock {\em Discrete Contin. Dyn. Syst. Ser. B}, 14(4):1465--1485, 2010.

\bibitem[Nat64]{Na64}
I.~P. Natanson.
\newblock {\em Constructive function theory. {V}ol. {I}. {U}niform
  approximation}.
\newblock Frederick Ungar Publishing Co., New York, 1964.
\newblock Translated from the Russian by Alexis N. Obolensky.

\bibitem[SSY25]{SSY25}
Pei Su, Chenmin Sun, and Xu~Yuan.
\newblock Quantitative observability for one-dimensional {S}chr\"odinger
  equations with potentials.
\newblock {\em J. Funct. Anal.}, 288(2):Paper No. 110695, 39, 2025.

\bibitem[Tao21]{Tao21}
Zhongkai Tao.
\newblock Exact control for the {S}chr\"odinger equation on the torus from
  small balls.
\newblock {\em Pure Appl. Anal.}, 3(2):387--401, 2021.

\bibitem[TT07]{TT07}
G{\'e}rald Tenenbaum and Marius Tucsnak.
\newblock New blow-up rates for fast controls of {Schr{\"o}dinger} and heat
  equations.
\newblock {\em Journal of Differential Equations}, 243(1):70--100, 2007.

\bibitem[TT09]{TT09}
G{\'e}rald Tenenbaum and Marius Tucsnak.
\newblock Fast and strongly localized observation for the {Schr{\"o}dinger}
  equation.
\newblock {\em Transactions of the American Mathematical Society},
  361(2):951--977, 2009.

\end{thebibliography}

\end{document}